\documentclass[12pt,letterpaper]{article}

\usepackage{renstyles}

\numberwithin{equation}{section}

\title{Reconstructing a state-independent cost function in a mean-field game model}
\author{
Kui Ren\thanks{Department Applied Physics and Applied Mathematics, Columbia University, New York, NY 10027; kr2002@columbia.edu}~\footnote{Author to whom any correspondence should be addressed}
\and
Nathan Soedjak\thanks{Department Applied Physics and Applied Mathematics, Columbia University, New York, NY 10027; ns3572@columbia.edu}
\and
Kewei Wang\thanks{School of Mathematical Sciences, Peking University, Beijing, China; wangkw@stu.pku.edu.cn}
\and
Hongyu Zhai\thanks{Department Applied Physics and Applied Mathematics, Columbia University, New York, NY 10027; hz2834@columbia.edu}
}

\begin{document}

\maketitle



\begin{abstract}
    In this short note, we consider an inverse problem to a mean-field games system where we are interested in reconstructing the state-independent running cost function from observed value-function data. We provide an elementary proof of a uniqueness result for the inverse problem using the standard multilinearization technique. One of the main features of our work is that we insist that the population distribution be a probability measure, a requirement that is not enforced in some of the existing literature on theoretical inverse mean-field games.
\end{abstract}





\section{Introduction}
\label{SEC:Intro}

This paper studies inverse problems to a mean-field game (MFG) model. Let $\Omega\subseteq\bbR^d$ ($d\ge 1$) be a bounded domain with smooth boundary $\partial\Omega$. Let $\bx\in\bbR^d$ denote the state variable, and $t\in[0,\infty)$ denote the time variable. We consider the mean-field game model in the following form
\begin{equation}\label{EQ:MFG}
    \left\{\begin{array}{lc}
        -\partial_t u-\Delta u + \frac{1}{2}|\nabla u|^2 = F(m(\bx,t)), &(\bx,t)\in \Omega \times (0,T)\\
        \partial_t m - \Delta m - \nabla\cdot(m\nabla u) = 0, &(\bx,t)\in \Omega \times (0,T)\\
        \partial_\nu u(\bx,t) = \partial_\nu m(\bx,t) = 0, &\bx\in\partial\Omega\\
        u(\bx,T) = 0,\ m(\bx,0)=m_0(\bx), &\bx\in \Omega
    \end{array}\right.
\end{equation}
where $u$ is the value function, $m$ is the population density and $F$ is the running cost function. In this work, we focus on the case when the Hamiltonian is quadratic. The terminal cost has been set to $0$, although it would be easy to add it back to the model by setting $u(\bx, T)=G(\bx)$ for some terminal cost function $G$. Without loss of generality, we let $\Omega$ have Lebesgue measure $1$.

The main feature of this model is that the running cost $F$ is independent of the state variable $\bx$ (as well as the time variable $t$), that is, $F=F(m)$. Such models of the cost functions can be viewed as simplified versions of the separable cost functions that can be written as a summation of a function of the state and a function of the population density; see, for instance, ~\cite{BaFe-NoDEA16,BaPr-SIAM14,GoPiSa-ESAIM16,LaWo-TRB-2011,PiVo-IUMJ17} for examples of such cost function models.

We are interested in reconstructing the cost function $F$ from the measurement of the total cost of the game. More precisely, we assume we have the data encoded in the map:
\begin{equation}\label{EQ:Data}
	\cM_{F}: m_0 \mapsto (V,W)
\end{equation}
 where $V$ and $W$ are the final cost of the game averaged over, respectively, the starting state and the players:
\begin{equation}\label{EQ:Data2}
V := \int_\Omega u(\bx,0)\,d\bx \quad \mbox{and}\quad W := \int_\Omega u(\bx,0)m_0(\bx)\,d\bx
\end{equation}
Our main purpose is to show that such data are sufficient to determine the running cost function $F(m)$ uniquely under mild additional assumptions. 

In recent years, there has been significant interest in studying inverse problems to mean-field game models. In most cases, the goal of such inverse problems is either to discover new forms of the model (for instance, to discover the running cost function, the terminal cost function, and the Hamiltonian) from observations or to calibrate known forms of the model (for instance, to determine parameters in given forms of the running cost function or the Hamiltonian) against measured data; see~\cite{AgLeFuNu-JCP22, ChFuLiNuOs-IP23,DiLiOsYi-JSC22,DiLiZh-arXiv23,FuLiOsLi-JCP23,ImLiYa-AML23,ImYa-arXiv23,KlAv-SIAM24,KlLiLi-arXiv23C,KlLiLi-JIIP24,KlLiYa-arXiv23,LeLiTeLiOs-SIAM21,LiMoZh-IP23,LiZh-arXiv23,LiJaLiNuOs-SIAM21,ReSoWa-arXiv23,YuXiChLa-arXiv24} and references therein for some recent analytical and numerical results in the area. 

Following previous works~\cite{LiMoZh-IP23,LiZh-arXiv22}, we restrict our interest to reconstructing $F$ in a special class of functions defined as follows.
\begin{definition}\label{DEF:Admissible}
    We say a function $U(z):\bbC\to\bbC$ is \emph{admissible} or in class $\cA$, that is, $U\in\cA$, if it satisfies the following conditions:\\[1ex]
    (i) The map $z\mapsto U(z)$ is holomorphic;\\[1ex]
    (ii) $U'(1)>0$.

    Therefore, if $U\in\cA$, then we have the power series representation 
    \begin{equation*}
        U(z) = \sum_{k=0}^{\infty} U^{(k)}\frac{(z-1)^k}{k!},\qquad z\in \bbR,
    \end{equation*}
    where $U^{(k)} = \frac{d^k U}{dz^k}(1)$ and $U^{(1)} > 0$.
\end{definition}
Let us mention that an admissible function $F\in \cA$ in the above definition is not required to satisfy $F(1)=0$. This is a critical difference from previous works based on the same expansion. 

For running cost functions in the admissible class, we can show the following local well-posedness result for the MFG system~\eqref{EQ:MFG}. For a given $F\in\cA$, let us define
\[
    u_b(\bx,t) := (T-t)F(1),\ \ \mbox{and}\ \  m_b(\bx,t) := 1\,.
\]
We observe that $(u_b, m_b)$ is a solution to~\eqref{EQ:MFG} with $m_0(\bx)=1$. Then~\eqref{EQ:MFG} is well-posed around $(u_b, m_b)$ in the following sense.
\begin{theorem}\label{THM:holomorphic}
    Let $Q:=\Omega\times[0,T]$ be the time-space domain. For a given $F\in \cA$, the following results hold:\\[1ex]
    (a) There exist constants $\delta,C>0$ such that for any
    \begin{equation}\nonumber
        m_0\in B_{1,\delta}\left(C^{2+\alpha}(\Omega)\right):=\{m_0\in C^{2+\alpha}(\Omega):\|m_0 - m_b(\cdot,0)\|_{C^{2+\alpha}(\Omega)}\le\delta\},
    \end{equation}
    the MFG system \eqref{EQ:MFG} admits a unique solution $(u,m)$ within the class
    \begin{equation}\nonumber
        \{(u,m)\in C^{2+\alpha,1+\frac{\alpha}{2}}(Q):\|(u,m) - (u_b, m_b)\|_{C^{2+\alpha,1+\frac{\alpha}{2}}(Q)}\le C\delta\}.
    \end{equation}
    Furthermore, the following estimate holds:
    \begin{equation}\nonumber
        \|(u,m) - (u_b, m_b)\|_{C^{2+\alpha,1+\frac{\alpha}{2}}(Q)}\le C\|m_0 - m_b(\cdot,0)\|_{C^{2+\alpha}(\Omega)}.
    \end{equation}

    (b) Consider the solution map
    \begin{equation}\nonumber
        S:B_{1,\delta}\left(C^{2+\alpha}(\Omega)\right)\to C^{2+\alpha,1+\frac{\alpha}{2}}(Q)
    \end{equation}
    defined by $S(m_0)=(u,m)$ where $(u,m)$ is the unique solution to system \eqref{EQ:MFG} in (a). Then $S$ is holomorphic.
\end{theorem}
This local well-posedness result can be proved in the same manner as the proof of~\cite[Theorem 3.1]{LiZh-arXiv22} using the implicit function theorem for Banach spaces. We will not reproduce the proof here. The most important point about this result is that the solution to the system is holomorphic with respect to $m_0(\bx)$ near $m_0(\bx) = 1$. This allows us to rigorously justify the linearization procedure we will use in Section~\ref{SEC:Proof} to prove the uniqueness result in Theorem~\ref{THM:Uniqueness}.

The rest of the paper is organized as follows. In Section~\ref{SEC:Uniqueness}, we state the uniqueness result of our inverse problem and highlight the main features that make our result different from existing ones in the literature. The proof of this result is presented in Section~\ref{SEC:Proof}. Concluding remarks are offered in Section~\ref{SEC:Conclusion}. 

\section{Uniqueness of reconstruction}
\label{SEC:Uniqueness}

We now state the uniqueness result we have for the inverse problem of determining $F(m)$ from data encoded in the map~\eqref{EQ:Data}. We restrict the domain of the map $\cM_{F}$ to lie in the set $\cP$ of probability density functions:
\begin{equation*}
    \cP := \left\{m: \Omega\to [0,\infty)\,:\, \int_\Omega m\,d\bx = 1\right\}.
\end{equation*}

The result is as follows.
\begin{theorem}\label{THM:Uniqueness}
Let $\Omega\subseteq\bbR^d$ ($d\ge 1$) be a bounded domain with smooth boundary $\partial\Omega$, and let $F_j\in \cA$ for $j=1,2$. Let $\cM_{F_j}$, $j=1,2$ be the measurement maps associated to the following system: 
\begin{equation}\label{EQ:MFG j}
    \left\{\begin{array}{lc}
        -\partial_t u_j - \Delta u_j +\frac12|\nabla u_j|^2 = F_j(m_j(\bx,t)), &(\bx,t)\in \Omega \times (0,T)\\
        \partial_t m_j - \Delta m_j - \nabla\cdot(m_j \nabla u_j) = 0, &(\bx,t)\in \Omega \times (0,T)\\
        \partial_\nu u_j(\bx,t) = \partial_\nu m_j(\bx,t) = 0, &\bx\in\partial\Omega\\
        u_j(\bx,T) = 0,\ m_j(\bx,0)=m_0(\bx), &\bx\in \Omega.
    \end{array}\right.
\end{equation}
If there exists $\delta>0$ such that for all $m_0\in \cP\cap B_{1,\delta}\left(C^{2+\alpha}(\Omega)\right)$ we have
\[
    \cM_{F_1}(m_0) = \cM_{F_2}(m_0)\ \ \Big(\mbox{i.e.}\  \left(V_1, W_1\right) = \left(V_2, W_2\right)\Big), 
\]
then it holds that 
\[
F_1(m) = F_2(m)\qquad\forall\ m\in \bbR.
\]
\end{theorem}

While this uniqueness result does not look particularly special, it does have a few features that make it different from existing results in the literature. First, in most of the existing results on the reconstruction of the running cost function, $F$ depends more than just the population density. One, therefore, often requires much more data than what we need here to have uniqueness for the inverse problem. Very often, data of the forms
\begin{equation*}
    \cM_F: m_0 \mapsto u(\bx, 0)\qquad \mbox{or} \qquad \cM_F: m_0 \mapsto \left(u(\bx, 0), m(\bx, T)\right)
\end{equation*}
are needed~\cite{LiMoZh-IP23,KlLiLi-arXiv23C}. Here, in our case, we need only two numbers $V$ and $W$ for each initial condition $m_0$ to uniquely determine the state-independent function $F(m)$. 

Secondly, since $m_0$ is a probability density function, it needs to be non-negative and normalized to $1$, that is, $m_0\in\cP$. This constraint, however, is not enforced in some of the existing works; see, for instance, ~\cite{LiMoZh-IP23, ReSoWa-arXiv23}. Even though initial conditions that are not probability density functions may be interpreted in some special settings, such as in the setting where $\Omega$ is just one of several disjoint domains on which the mean-field game is played, $m_0$ being a probability density function is more natural for many applications. Enforcing that $m_0\in\cP$, however, makes the inverse problem significantly more challenging. For instance, to continue using the multilinearization techniques in~\cite{LiMoZh-IP23, ReSoWa-arXiv23}, we need to be able to find an explicit background solution $(u_b, m_b)$ around which we can linearize the forward model. Such a background solution is not easy to find unless we know the evaluation of the running cost function (which we do not know). This is why $F$ is assumed to vanish when $m_b=1$ in~\cite{LiZh-arXiv22}.  In our case here, $F(1)$ is allowed to be nonzero.

It is worth noting that the boundary conditions $\partial_\nu u = \partial_\nu m = 0$ on $\partial\Omega$ represent reflection at the boundary $\partial\Omega$, and thus the total population is conserved over time. This can be seen with an explicit calculation using integration by parts:
\begin{align*}
    \frac{d}{dt}\int_{\Omega} m\,d\bx &= \int_{\Omega} \Delta m + \nabla\cdot(m\nabla u)\, d\bx \\
    &=\int_{\partial\Omega} (\partial_\nu m + m\partial_\nu u)\, dS = 0. 
\end{align*}

\section{Proof of \texorpdfstring{Theorem~\ref{THM:Uniqueness}}{}}
\label{SEC:Proof}

We now provide an elementary proof of the uniqueness result. Our proof is based on the standard multilinearization technique adopted in~\cite{LiMoZh-IP23,LiZh-arXiv22,ReSoWa-arXiv23}. 

Without loss of generality, we let $\Omega$ have Lebesgue measure $1$. For $\ell=1,\dots,N$, let $g_\ell\in C^{2+\alpha}(\Omega)$ be an arbitrary function with $\int_\Omega g_\ell = 0$. Set 
\begin{equation*}
m_0(\bx) = 1+\sum_{\ell=1}^N \eps_\ell g_\ell(\bx), 
\end{equation*}
and observe that for sufficiently small $\eps_\ell$, the probability density constraint on $m_0$ is satisfied. Now note that when $\eps := (\eps_1,\dots,\eps_N) = 0$, we have $m_j(\bx,t) = 1$ and $u_j(\bx,t) = (T-t)F_j(1)$. In particular, $\nabla u_j = 0$ when $\eps = 0$. By plugging in $t=0$ and integrating over $\Omega$, we can recover $F(1)$ from the data $\cM(1)$: 
\begin{align*}
    F_1(1) = \frac{1}{T}\int_{\Omega} u_1(\bx, 0)\,d\bx = \frac{1}{T}\int_{\Omega} u_2(\bx,0)\,d\bx = F_2(1).
\end{align*}

Define $u_j^{(\ell)}(\bx,t) := \partial_{\eps_\ell} u_j(\bx,t;\eps) |_{\eps = 0}$, $m_j^{(\ell)}(\bx,t) := \partial_{\eps_\ell} m_j(\bx,t;\eps) |_{\eps = 0}$, $F_j^{(1)} := \partial_m F_j(m)|_{m=1}$. More generally, define
\begin{equation*}
    \begin{aligned}
        u_j^{(\ell_1,\ell_2,\dots,\ell_M)}(\bx,t) &:= \partial_{\varepsilon_{\ell_1}}\partial_{\varepsilon_{\ell_2}}\cdots\partial_{\varepsilon_{\ell_M}}u_j(\bx,t;\varepsilon)|_{\varepsilon=0},\\
        m_j^{(\ell_1,\ell_2,\dots,\ell_M)}(\bx,t) &:= \partial_{\varepsilon_{\ell_1}}\partial_{\varepsilon_{\ell_2}}\cdots\partial_{\varepsilon_{\ell_M}}m_j(\bx,t;\varepsilon)|_{\varepsilon=0},\\
        F_j^{(M)} &:= \partial_m^M F_j(m)|_{m=1}.
    \end{aligned}
\end{equation*}
These derivatives exist by the infinite differentiability of the solution map $S$ established in \Cref{THM:holomorphic}.

\subsection{First-order linearization}

By differentiating the MFG system \eqref{EQ:MFG j} with respect to $\eps_\ell$ at $\eps = 0$, we find that
\begin{equation}\label{EQ:Ord1}
    \left\{\begin{array}{lc}
        -\partial_t u_j^{(\ell)} - \Delta u_j^{(\ell)} = F_j^{(1)} m_j^{(\ell)}, &(\bx,t)\in \Omega \times (0,T),\\
        \partial_t m_j^{(\ell)} - \Delta m_j^{(\ell)} - \Delta u_j^{(\ell)} = 0, &(\bx,t)\in \Omega \times (0,T),\\
        \partial_\nu u_j^{(\ell)}(\bx,t) = \partial_\nu m_j^{(\ell)}(\bx,t) = 0, &\bx\in\partial\Omega,\\
        u_j^{(\ell)}(\bx,T) = 0,\ m_j^{(\ell)}(\bx,0) = g_\ell(\bx), &\bx\in \Omega.
    \end{array}\right.
\end{equation}
Observe that $F_j^{(1)} > 0$ due to the admissibility assumption (recall \Cref{DEF:Admissible}).

To start, consider the two differential data measurements 
\[
W_j^{(\ell,\ell)} := \partial_{\varepsilon_{\ell}}^2 W_j(\varepsilon)|_{\varepsilon=0},\qquad V_j^{(\ell,\ell)} := \partial_{\varepsilon_{\ell}}^2 V_j(\varepsilon)|_{\varepsilon=0}.
\]
We compute
\begin{align*}
    W_j^{(\ell,\ell)} &= \int_{\Omega} u_j^{(\ell,\ell)}(\bx,0)\,d\bx + 2\int_{\Omega} u_j^{(\ell)}(\bx,0)m_j^{(\ell)}(\bx,0)\,d\bx + \int_{\Omega}u_j(\bx,0)m_j^{(\ell,\ell)}(\bx,0)\,d\bx\\
    &= V_j^{(\ell,\ell)} + 2\int_{\Omega} u_j^{(\ell)}(\bx,0)m_j^{(\ell)}(\bx,0)\,d\bx,
\end{align*}
using the fact that $m_j^{}(\bx, 0)=0$. Therefore, the quantity
\begin{align}\label{EQ:Data Ord1}
    \int_{\Omega} u_j^{(\ell)}(\bx,0)m_j^{(\ell)}(\bx,0)\,d\bx = \frac{W_j^{(\ell,\ell)} - V_j^{(\ell,\ell)}}{2}
\end{align}
is independent of $j\in\{1,2\}$ from the data.

Next, let $(\phi_n(\bx))_{n=0}^{\infty}$ be the $L^2$-normalized Neumann eigenfunctions of $-\Delta$ on $\Omega$, with corresponding eigenvalues $0 = \lambda_0 < \lambda_1 \le \lambda_2 \le \cdots $. We now perform eigenfunction expansion in this basis to convert the above system of PDEs into a system of ODEs. Let
\begin{align*}
    u_j^{(\ell)}(\bx,t) = \sum_{n=0}^{\infty} U_n^j(t) \phi_n(\bx), \qquad m_j^{(\ell)}(\bx,t) = \sum_{n=0}^{\infty} M_n^j(t) \phi_n(\bx).
\end{align*}
Since the eigenfunctions are orthogonal (and $1 = \phi_0(\bx)$ is an eigenfunction), the constraint that $\int_\Omega g_\ell = 0$ is then equivalent to 
\begin{align*}
M_0^j(0) = 0.
\end{align*}
From the zero terminal condition on $u_j^{(\ell)}$ we also see that for all $n\ge 0$, 
\begin{align*}
U_n^j(T) = 0. 
\end{align*}

Now, by differentiating term-by-term, we see that the PDE system turns into the ODE system, for each $n\ge 0$,
\begin{align*}
    -\frac{d U_n^j}{dt} + \lambda_n U_n^j &= F_j^{(1)} M_n^j,\\
    \frac{d M_n^j}{dt} + \lambda_n M_n^j + \lambda_n U_n^j &= 0.
\end{align*}
For $n=0$, we have $\lambda_0 = 0$, and then the boundary conditions imply that $M_0^j(t) = U_0^j(t) = 0$ for all $t$. 

For $n > 0$, we write the ODE system in matrix form  
\begin{align}\label{EQ:ODE}
    \frac{d}{dt} \begin{pmatrix}U_n^j\\ M_n^j\end{pmatrix} = A \begin{pmatrix}U_n^j\\ M_n^j\end{pmatrix},
\end{align}
where 
\begin{align*}
    A:= \begin{pmatrix}\lambda_n & -F_j^{(1)}\\ 
    -\lambda_n & -\lambda_n\end{pmatrix} 
    = S
    \begin{pmatrix}\sqrt{\lambda_n(\lambda_n + F_j^{(1)})} & 0 \\ 
    0 & -\sqrt{\lambda_n(\lambda_n+F_j^{(1)})}\end{pmatrix} 
    S^{-1},
\end{align*}
where 
\begin{align*}
    S &= \begin{pmatrix}-\sqrt{\lambda_n}-\sqrt{\lambda_n + F_j^{(1)}} & -\sqrt{\lambda_n} + \sqrt{\lambda_n + F_j^{(1)}} \\ 
    \sqrt{\lambda_n} & \sqrt{\lambda_n}\end{pmatrix},\\
    S^{-1} &= \frac{1}{2\sqrt{\lambda_n(\lambda_n + F_j^{(1)})}}
    \begin{pmatrix}-\sqrt{\lambda_n} & -\sqrt{\lambda_n} + \sqrt{\lambda_n + F_j^{(1)}} \\ 
    \sqrt{\lambda_n} & \sqrt{\lambda_n} + \sqrt{\lambda_n + F_j^{(1)}}\end{pmatrix}.
\end{align*}
With this diagonalization in hand, we can explicitly solve the ODE:
\begin{align*}
    \begin{pmatrix}U_n^j(T) \\ M_n^j(T)\end{pmatrix} = S
    \begin{pmatrix}\exp\left(\sqrt{\lambda_n(\lambda_n + F_j^{(1)})} T\right)& 0 \\ 
    0 & \exp\left(-\sqrt{\lambda_n(\lambda_n+F_j^{(1)})} T\right)\end{pmatrix} 
    S^{-1} \begin{pmatrix}U_n^j(0) \\ M_n^j(0)\end{pmatrix},
\end{align*}
or in reverse:
\begin{align*}
    \begin{pmatrix}U_n^j(0) \\ M_n^j(0)\end{pmatrix} = S
    \begin{pmatrix}\exp\left(-\sqrt{\lambda_n(\lambda_n + F_j^{(1)})} T\right)& 0 \\ 
    0 & \exp\left(\sqrt{\lambda_n(\lambda_n+F_j^{(1)})} T\right)\end{pmatrix} 
    S^{-1} \begin{pmatrix}0 \\ M_n^j(T)\end{pmatrix}.
\end{align*}
(Recall that $U_n^j(T) = 0$.) Multiplying out the right hand side leads to
\begin{align*}
\begin{pmatrix}U_n^j(0) \\ M_n^j(0)\end{pmatrix} = \frac{M_n^j(T)}{2\sqrt{\lambda_n(\lambda_n + F_j^{(1)})}}
\begin{pmatrix}
-F_j^{(1)}\exp\left(-\sqrt{\lambda_n(\lambda_n + F_j^{(1)})} T\right) + F_j^{(1)}\exp\left(\sqrt{\lambda_n(\lambda_n + F_j^{(1)})} T\right)\\
A_{n,j}\exp\left(-\sqrt{\lambda_n(\lambda_n + F_j^{(1)})} T\right) + B_{n,j}\exp\left(\sqrt{\lambda_n(\lambda_n + F_j^{(1)})} T\right)
\end{pmatrix},
\end{align*}
where 
\begin{align*}
    A_{n,j} &:= -\lambda_n + \sqrt{\lambda_n(\lambda_n + F_j^{(1)})} \ge 0,\\
    B_{n,j} &:= \lambda_n + \sqrt{\lambda_n(\lambda_n + F_j^{(1)})} > 0.
\end{align*}
(More generally, we have that the full solution $\begin{pmatrix}U_n^j(t) \\ M_n^j(t)\end{pmatrix}$ to the ODE \eqref{EQ:ODE} is
\begin{align}\label{EQ:ODE Sol}
\frac{M_n^j(T)}{2\sqrt{\lambda_n(\lambda_n + F_j^{(1)})}}
\begin{pmatrix}
-F_j^{(1)}\exp\left(-\sqrt{\lambda_n(\lambda_n + F_j^{(1)})} (T-t)\right) + F_j^{(1)}\exp\left(\sqrt{\lambda_n(\lambda_n + F_j^{(1)})} (T-t)\right)\\
A_{n,j}\exp\left(-\sqrt{\lambda_n(\lambda_n + F_j^{(1)})} (T-t)\right) + B_{n,j}\exp\left(\sqrt{\lambda_n(\lambda_n + F_j^{(1)})} (T-t)\right)
\end{pmatrix}.
\end{align}
In particular, either (i) $M_n^j(t) = 0$ for all $t$,  or (ii) $M_n^j(t)$ is strictly positive for all $t$ or strictly negative for all $t$. We will use this fact later in the second-order and higher-order linearization steps.)

If $M_n^j(T) = 0$, we have $U_n^j(0) = M_n^j(0) = 0$, so $U_n^j(t) = M_n^j(t) = 0$ for all $t$. On the other hand, if $M_n^j(T) \neq 0$, then $M_n^j(0) \neq 0$ and 
\begin{align*}
    \frac{U_n^j(0)}{M_n^j(0)} = \frac{-F_j^{(1)} + F_j^{(1)}\exp\left(2\sqrt{\lambda_n(\lambda_n + F_j^{(1)})} T\right)}{A_{n,j} + B_{n,j}\exp\left(2\sqrt{\lambda_n(\lambda_n + F_j^{(1)})} T\right)}.
\end{align*}
It follows that for any $a\in \bbR$, there is a unique solution to the ODE \eqref{EQ:ODE} with the initial and terminal boundary conditions 
\begin{align*}
    M_n^j(0) = a,\qquad U_n^j(T) = 0,
\end{align*}
and this solution has the initial data
\begin{align*}
    U_n^j(0) = a\cdot \frac{-F_j^{(1)} + F_j^{(1)}\exp\left(2\sqrt{\lambda_n(\lambda_n + F_j^{(1)})} T\right)}{A_{n,j} + B_{n,j}\exp\left(2\sqrt{\lambda_n(\lambda_n + F_j^{(1)})} T\right)}.
\end{align*}
Now fix an integer $k\ge 1$ (e.g., $k=1$ will do), and consider the following initial condition for $m_j^{(\ell)}$:
\begin{align*}
    m_j^{(\ell)}(\bx, 0) = g_\ell(\bx) = \phi_k(\bx).
\end{align*}
(Note that $\int_\Omega g_\ell = 0$.) That is, we set $M_n^j(0)$ to be $1$ if $n=k$, and $0$ otherwise. Then, by what we just showed, the unique solution to \eqref{EQ:Ord1} has initial data $U_n^j(0)$ equal to 
\begin{align*}
    \frac{-F_j^{(1)} + F_j^{(1)}\exp\left(2\sqrt{\lambda_n(\lambda_n + F_j^{(1)})} T\right)}{A_{n,j} + B_{n,j}\exp\left(2\sqrt{\lambda_n(\lambda_n + F_j^{(1)})} T\right)}
\end{align*}
if $n=k$, and $0$ otherwise. But recall from the data \eqref{EQ:Data Ord1} that
\begin{align*}
    \int_{\Omega} u_j^{(\ell)}(\bx,0)m_j^{(\ell)}(\bx,0)\,d\bx = \frac{W^{(\ell,\ell)} - V^{(\ell,\ell)}}{2}
\end{align*}
is independent of $j$, which in this case means 
\begin{align*}
    U_k^j(0) = \int_{\Omega} u_j^{(\ell)}(\bx,0)\phi_k(\bx)\,d\bx = \frac{W^{(\ell,\ell)} - V^{(\ell,\ell)}}{2}
\end{align*}
is independent of $j$. That is, 
\begin{align*}
    \frac{-F_j^{(1)} + F_j^{(1)}\exp\left(2\sqrt{\lambda_k(\lambda_k + F_j^{(1)}}) T\right)}{\left(-\lambda_k + \sqrt{\lambda_k(\lambda_k + F_j^{(1)})}\right) + \left(\lambda_k + \sqrt{\lambda_k(\lambda_k + F_j^{(1)})}\right)\exp\left(2\sqrt{\lambda_k(\lambda_k + F_j^{(1)})} T\right)}
\end{align*}
is independent of $j\in \{1,2\}$. So to prove that $F_1^{(1)} = F_2^{(1)}$, it suffices to show that the above expression is a strictly increasing function of $F_j^{(1)}$ on $(0, \infty)$. The proof goes as follows.

Employing the change of variable
\begin{align*}
G_{k}=\sqrt{\lambda_{k} (\lambda_{k} +F^{(1)}_{j} )}, 
\end{align*}
the expression above can be written as
\begin{align*}
\frac{\frac{G^{2}_{k}-\lambda^{2}_{k} }{\lambda_{k} } \left( \exp \left( 2TG_{k}\right)  -1\right)  }{\left( G_{k}-\lambda_{k} \right)  +\left( \lambda_{k} +G_{k}\right)  \exp \left( 2TG_{k}\right)  }.
\end{align*}
Its derivative with respect to $G_k$ is
\begin{align*}
&\frac{d}{dG_{k}} \frac{\frac{G^{2}_{k}-\lambda^{2}_{k} }{\lambda_{k} } \left( \exp \left( 2TG_{k}\right)  -1\right)  }{\left( G_{k}-\lambda_{k} \right)  +\left( \lambda_{k} +G_{k}\right)  \exp \left( 2TG_{k}\right)  } \\
=&\frac{\left( \lambda_{k} +G_{k}\right)^{2}  e^{4TG_{k}}+4TG_{k}\left( G^{2}_{k}-\lambda^{2}_{k} \right)  e^{2TG_{k}}-4\lambda_{k} G_{k}e^{2TG_{k}}-\left( G_{k}-\lambda_{k} \right)^{2}  }{\lambda_{k} \left( \left( G_{k}-\lambda_{k} \right)  +\left( \lambda_{k} +G_{k}\right)  e^{2TG_{k}}\right)^{2}  }. 
\end{align*}
Notice that the denominator is strictly positive given $\lambda_{k} >0$, which means we can rewrite the numerator as
\begin{align*}
&\left( \lambda_{k} +G_{k}\right)^{2}  e^{4TG_{k}}+4TG_{k}\left( G^{2}_{k}-\lambda^{2}_{k} \right)  e^{2TG_{k}}-4\lambda_{k} G_{k}e^{2TG_{k}}-\left( G_{k}-\lambda_{k} \right)^{2}  \\
=&\left( \lambda_{k} +G_{k}\right)^{2}  \left( e^{4TG_{k}}-1\right)  -4\lambda_{k} G_{k}\left( e^{2TG_{k}}-1\right)  +4TG_{k}\left( G^{2}_{k}-\lambda^{2}_{k} \right)  e^{2TG_{k}}\\
=&\left( e^{2TG_{k}}-1\right)  \left( \left( \lambda_{k} +G_{k}\right)^{2}  \left( e^{2TG_{k}}+1\right)  -4\lambda_{k} G_{k}\right)  +4TG_{k}\left( G^{2}_{k}-\lambda^{2}_{k} \right)  e^{2TG_{k}}\\
\ge &\left( e^{2TG_{k}}-1\right)  \left( 2\left( \lambda_{k} +G_{k}\right)^{2}  -4\lambda_{k} G_{k}\right)  +4TG_{k}\left( G^{2}_{k}-\lambda^{2}_{k} \right)  e^{2TG_{k}}\\
\ge&\ 0 .
\end{align*}
On the other hand, it holds that
\begin{align*}
\frac{dG_{k}}{dF^{(1)}_{j} } &=\frac{d}{dF^{(1)}_{j} } \sqrt{\lambda_{k} (\lambda_{k} +F^{(1)}_{j} )} =\frac{\lambda_{k} }{2\sqrt{\lambda_{k} \left( \lambda_{k} +F^{(1)}_{j}\right)  } } >0 .
\end{align*}
Hence, we conclude that
\begin{align*}
\frac{d}{dF^{(1)}_{j} } \frac{F^{(1)}_{j} \left( \exp \left( 2T\sqrt{\lambda_{k} (\lambda_{k} + F^{(1)}_{j} )} \right)  -1\right)  }{\left( \sqrt{\lambda_{k} (\lambda_{k} +F^{(1)}_{j} )} -\lambda_{k} \right)  +\left( \lambda_{k} +\sqrt{\lambda_{k} (\lambda_{k} + F^{(1)}_{j} )} \right)  \exp \left( 2T\sqrt{\lambda_{k} (\lambda_{k} + F^{(1)}_{j} )} \right)  } > 0.
\end{align*}
We have, therefore, shown that 
\begin{equation*}
F_1^{(1)} = F_2^{(1)}.
\end{equation*}
It follows from this that $u_1^{(\ell)} = u_2^{(\ell)}$ and $m_1^{(\ell)} = m_2^{(\ell)}$ as well. This is because the above arguments show that the system \eqref{EQ:Ord1} has a unique solution.

\subsection{Second-order linearization}

By differentiating the MFG system \eqref{EQ:MFG j} with respect to $\eps_1$ and $\eps_2$ at $\eps = 0$, we find that: 

\begin{equation}\label{EQ:Ord2}
    \left\{\begin{array}{lc}
        -\partial_t u_j^{(1,2)} - \Delta u_j^{(1,2)} + \nabla u_j^{(1)}\cdot \nabla u_j^{(2)} = F_j^{(2)} m_j^{(1)} m_j^{(2)} + F_j^{(1)}m_j^{(1,2)}, &(\bx,t)\in \Omega \times (0,T),\\
        \partial_t m_j^{(1,2)} - \Delta m_j^{(1,2)} - \Delta u_j^{(1,2)} = \nabla\cdot(m_j^{(1)}\nabla u_j^{(2)}) + \nabla\cdot (m_j^{(2)}\nabla u_j^{(1)}), &(\bx,t)\in \Omega \times (0,T),\\
        \partial_\nu u_j^{(1,2)}(\bx,t) = \partial_\nu m_j^{(1,2)}(\bx,t) = 0, &\bx\in\partial\Omega,\\
        u_j^{(1,2)}(\bx,T) = 0,\ m_j^{(1,2)}(\bx,0) = 0, &\bx\in \Omega.
    \end{array}\right.
\end{equation}
Define $\overline{u} := u_1^{(1,2)} - u_2^{(1,2)}$, $\overline{m} := m_1^{(1,2)} - m_2^{(1,2)}$, $\overline{F} := F_1^{(2)} - F_2^{(2)}$. Subtracting the above system for $j=1,2$ and using the fact that the first-order terms do not depend on $j$ (thanks to the first-order linearization step) gives 
\begin{equation}\label{EQ:Ord2 bar}
    \left\{\begin{array}{lc}
        -\partial_t \overline{u} - \Delta \overline{u} = \overline{F} m^{(1)} m^{(2)} + F^{(1)}\overline{m}, &(\bx,t)\in \Omega \times (0,T),\\
        \partial_t \overline{m} - \Delta \overline{m} - \Delta \overline{u} = 0, &(\bx,t)\in \Omega \times (0,T),\\
        \partial_\nu \overline{u}(\bx,t) = \partial_\nu \overline{m}(\bx,t) = 0, &\bx\in\partial\Omega,\\
        \overline{u}(\bx,T) = 0,\ \overline{m}(\bx,0) = 0, &\bx\in \Omega.
    \end{array}\right.
\end{equation}
As before, we perform eigenfunction expansion. Let
\begin{align*}
    \overline{u}(\bx,t) = \sum_{n=0}^{\infty} \overline{U}_n(t) \phi_n(\bx), \qquad \overline{m}(\bx,t) = \sum_{n=0}^{\infty} \overline{M}_n(t) \phi_n(\bx).
\end{align*}
Then for each $n\ge 0$, we obtain the ODE system
\begin{align*}
    -\frac{d \overline{U}_n}{dt} + \lambda_n \overline{U}_n &= F^{(1)} \overline{M}_n + \overline{F}\langle m^{(1)}m^{(2)}, \phi_n\rangle_{L^2(\Omega)},\\
    \frac{d \overline{M}_n}{dt} + \lambda_n \overline{M}_n + \lambda_n \overline{U}_n &= 0.
\end{align*}
We will, in fact, only need to consider this system for $n=0$. In this case $\lambda_n = 0$, so the system simplifies to
\begin{align*}
    -\frac{d \overline{U}_0}{dt} &= F^{(1)} \overline{M}_0 + \overline{F}\langle m^{(1)}m^{(2)}, 1\rangle_{L^2(\Omega)},\\
    \frac{d \overline{M}_0}{dt} &= 0.
\end{align*}
Combined with the initial condition $\overline{M}_0(0) = 0$ (from $\overline{m}(\bx,0) = 0$), this gives $\overline{M}_0(t) = 0$ for all $t$. 

Now note that from the data, 
\begin{align*}
    \int_{\Omega} u_j^{(1,2)}(\bx,0)\,d\bx = V^{(1,2)}
\end{align*}
is independent of $j$. In other words, 
\begin{align*}
    \overline{U}_0(0) = \int_{\Omega} \overline{u}(\bx,0)\cdot 1\,d\bx = 0.
\end{align*}
Also, from the terminal condition $\overline{u}(\bx, T) = 0$ we have $\overline{U}_0(T) = 0$. Therefore, 
\begin{align}\label{EQ:Nonzero Int}
0 = \overline{U}_0(T) - \overline{U}_0(0) = \int_0^T \frac{d \overline{U}_0}{dt} \,dt = -\overline{F} \int_0^T \langle m^{(1)}m^{(2)}, 1\rangle_{L^2(\Omega)} \,dt. 
\end{align}
We want to make the integral nonzero so that we can conclude that $\overline{F} = 0$.

Let us set
\begin{align*}
    g_1(\bx) = g_2(\bx) = \phi_1(\bx).
\end{align*}
(This is allowed because as a non-constant eigenfunction, $\phi_1(\bx)$ satisfies $\int_\Omega \phi_1(\bx)\,d\bx = 0$.) Then, by the discussion surrounding \eqref{EQ:ODE Sol} in the first-order linearization step, we have that 
\begin{align*}
    m^{(1)}(\bx,t) = m^{(2)}(\bx,t) = M_1(t)\phi_1(\bx),
\end{align*}
for some strictly positive function $M_1(t) > 0$.  The integral in \eqref{EQ:Nonzero Int} is therefore strictly positive: 
\begin{align*}
    \int_0^T \langle m^{(1)}m^{(2)}, 1\rangle_{L^2(\Omega)}\,dt = \int_0^T M_1^2(t)\,dt \int_\Omega \phi_1^2(\bx)\,d\bx > 0. 
\end{align*}
We conclude from \eqref{EQ:Nonzero Int} that $\overline{F} = 0$. That is, we have shown 
\begin{align*}
    F_1^{(2)} = F_2^{(2)}.
\end{align*}
It follows from this that $\overline{u}=\overline{m}=0$, i.e., $u_1^{(1,2)} = u_2^{(1,2)}$ and $m_1^{(1,2)} = m_2^{(1,2)}$ as well. This is because the ODE arguments in the first-order linearization step show that the system \eqref{EQ:Ord2 bar} has a unique solution.

\subsection{Higher-order linearization}
We first introduce an auxiliary lemma as follows.

\begin{lemma}\label{LEM:Eig}
    Let $(\phi_k(\bx))_{k=0}^{\infty}$ be the $L^2$-normalized Neumann eigenfunctions of $-\Delta$ on $\Omega$, with corresponding eigenvalues $0 = \lambda_0 < \lambda_1 \le \lambda_2 \le \cdots $. Then for all $n\ge 2$, there exist $k_1,\dots,k_n \ge 1$ such that $\langle \phi_{k_1}\phi_{k_2} \cdots \phi_{k_n}, 1\rangle_{L^2(\Omega)}\neq 0$.
\end{lemma}

\begin{proof}
    In fact, we can fix $k_1=k_2=\cdots=k_{n-1}=1$. Assume that for all $k\ge 1$, $\langle \phi_1^{n-1}\phi_{k}, 1\rangle_{L^2(\Omega)}=0$. Since $(\phi_k(\bx))_{k=1}^{\infty}$ is a complete orthonormal basis of
    \begin{equation*}
        L_0^2(\Omega):=\left\{f\in L^2(\Omega):\int_{\Omega}f d\bx=0\right\},
    \end{equation*}
    we have $\langle \phi_1^{n-1}f, 1\rangle_{L^2(\Omega)}=0$ for all $f\in L_0^2(\Omega)$. This implies $\phi_1^{n-1}$ is a constant function, which contradicts the fact that $\phi_1$ is non-constant. Therefore, there exists $k_n\ge 1$ such that $\langle \phi_1^{n-1}\phi_{k_n}, 1\rangle_{L^2(\Omega)}\neq 0$.
\end{proof}

By differentiating the MFG system \eqref{EQ:MFG j} with respect to $\eps_1, \dots, \eps_n$ at $\eps = 0$, we find that:

\begin{equation}\label{EQ:Ordn}
    \left\{\begin{array}{lc}
        -\partial_t u_j^{(1,\dots,n)} - \Delta u_j^{(1,\dots,n)} = F_j^{(n)} m_j^{(1)} m_j^{(2)} \cdots m_j^{(n)} + F_j^{(1)}m_j^{(1,\dots,n)} + l.o.t. &(\bx,t)\in \Omega \times (0,T)\\
        \partial_t m_j^{(1,\dots,n)} - \Delta m_j^{(1,\dots,n)} - \Delta u_j^{(1,\dots,n)} = l.o.t., &(\bx,t)\in \Omega \times (0,T)\\
        \partial_\nu u_j^{(1,\dots,n)}(\bx,t) = \partial_\nu m_j^{(1,\dots,n)}(\bx,t) = 0, &\bx\in\partial\Omega\\
        u_j^{(1,\dots,n)}(\bx,T) = 0,\ m_j^{(1,\dots,n)}(\bx,0) = 0, &\bx\in \Omega
    \end{array}\right.
\end{equation}
where ``l.o.t." stands for lower-order terms in the linearization process.

We proceed inductively. Let $n\ge 2$, and assume that all lower-order terms are independent of $j$. (In particular, this means $F_1^{(k)} = F_2^{(k)}$ for all $k\le n-1$.)

Define $\overline{u} := u_1^{(1,\dots,n)} - u_2^{(1,\dots,n)}$, $\overline{m} := m_1^{(1,\dots,n)} - m_2^{(1,\dots,n)}$, $\overline{F} := F_1^{(n)} - F_2^{(n)}$. Subtracting the above system for $j=1,2$ and using the fact that the lower-order terms do not depend on $j$ gives
\begin{equation}\label{EQ:Ordn bar}
    \left\{\begin{array}{lc}
        -\partial_t \overline{u} - \Delta \overline{u} = \overline{F} m^{(1)} m^{(2)}\cdots m^{(n)} + F^{(1)}\overline{m}, &(\bx,t)\in \Omega \times (0,T),\\
        \partial_t \overline{m} - \Delta \overline{m} - \Delta \overline{u} = 0, &(\bx,t)\in \Omega \times (0,T),\\
        \partial_\nu \overline{u}(\bx,t) = \partial_\nu \overline{m}(\bx,t) = 0, &\bx\in\partial\Omega,\\
        \overline{u}(\bx,T) = 0,\ \overline{m}(\bx,0) = 0, &\bx\in \Omega.
    \end{array}\right.
\end{equation}
Following exactly the same reasoning as in the second-order linearization step leads to the identity
 \begin{align}\label{EQ:Nonzero Int Ordn}
0 = \overline{U}_0(T) - \overline{U}_0(0) = \int_0^T \frac{d \overline{U}_0}{dt} \,dt = -\overline{F} \int_0^T \langle m^{(1)}m^{(2)}\cdots m^{(n)}, 1\rangle_{L^2(\Omega)} \,dt. 
\end{align}
In order to make the integral nonzero, we will choose $g_1,g_2,\dots,g_n$ in the following way. By \Cref{LEM:Eig}, there exist $k_1,\dots,k_n \ge 1$ such that $\langle \phi_{k_1}\phi_{k_2} \cdots \phi_{k_n}, 1\rangle_{L^2(\Omega)}\neq 0$. For $\ell=1,\dots, n$, we now choose 
\begin{equation*}
    g_\ell(\bx) = \phi_{k_\ell}(\bx).
\end{equation*}
(This is allowed because as a non-constant eigenfunction, $\phi_{k_\ell}(\bx)$ satisfies $\int_\Omega \phi_{k_\ell}(\bx)\,d\bx = 0$.) Then, by the discussion surrounding~\eqref{EQ:ODE Sol} in the first-order linearization step, we have that 
\begin{equation*}
m^{(\ell)}(\bx,t) = M_\ell(t)\phi_{k_\ell}(\bx),
\end{equation*}
for some strictly positive functions $M_\ell(t) > 0$ ($\ell=1,\dots,n$). Hence,
\begin{equation*}
    m^{(1)}(\bx,t) m^{(2)}(\bx,t)\cdots m^{(n)}(\bx,t) = M_1(t) M_2(t) \cdots M_n(t) \phi_{k_1}(\bx)\phi_{k_2}(\bx) \cdots \phi_{k_n}(\bx)
\end{equation*}
The integral in \eqref{EQ:Nonzero Int Ordn} is therefore nonzero: 
\begin{align*}
\int_0^T \langle m^{(1)}m^{(2)}\cdots m^{(n)}, 1\rangle_{L^2(\Omega)} \,dt = \langle \phi_{k_1}\phi_{k_2} \cdots \phi_{k_n}, 1\rangle_{L^2(\Omega)}\int_0^T M_1(t)M_2(t)\dots M_n(t)\,dt \neq 0.
\end{align*}
We conclude from \eqref{EQ:Nonzero Int Ordn} that $\overline{F} = 0$. That is, we have shown 
\begin{align*}
    F_1^{(n)} = F_2^{(n)}.
\end{align*}
It follows from this that $\overline{u}=\overline{m}=0$, i.e., $u_1^{(1,\dots,n)} = u_2^{(1,\dots,n)}$ and $m_1^{(1,\dots,n)} = m_2^{(1,\dots,n)}$ as well. This is because the ODE arguments in the first-order linearization step show that the system \eqref{EQ:Ordn bar} has a unique solution. 

The induction is thus complete. We have shown that 
\begin{align*}
    F_1^{(n)} = F_2^{(n)},\quad \forall\, n\ge 0.
\end{align*}
Hence, $F_1 = F_2$, and the proof is finished. \hfill $\qed$

\section{Concluding remarks}
\label{SEC:Conclusion}

We considered an inverse problem to a mean-field games system where we are interested in reconstructing a state-independent running cost function from observed value-function data. We provide an elementary proof of a uniqueness result for the inverse problem using the standard multilinearization technique. 

By assuming that the running cost is independent of the state variable, we made it possible to perform multilinearization around a special background state of $m_b=1$. Moreover, such a linearization around $m_b=1$ allows us to enforce the constraint that the initial condition $m_0$ is indeed a population density function, making the uniqueness theory meaningful in more general circumstances than those previously reported.

There are still many seemingly simple inverse problems to be solved for the mean-field game system~\eqref{EQ:MFG}. One such problem is to reconstruct a running cost function of the form $F(\bx, m)=G(m)+f(\bx)$ with $G$ and $f$ both unknown. While this problem seems like a small perturbation of the problem we solved, it is overwhelmingly challenging as the method of our proof fails completely due to the lack of a known background solution pair $(u_b, m_b)$ around which we can linearized the MFG system. Even the case of known $G(m)$, that is, the case where only $f(\bx)$ is unknown, is very interesting.

\section*{Acknowledgments}

We would like to thank the anonymous referees for their useful comments, especially for providing many references~\cite{ImLiYa-IPI24,KlAv-SIAM24,KlLiLi-JIIP24,LiZh-arXiv24}, that helped us improve the quality of this work. 

This work is partially supported by the National Science Foundation through grants DMS-1937254 and DMS-2309802.
    
{\small

\begin{thebibliography}{10}

\bibitem{AgLeFuNu-JCP22}
{\sc S.~Agrawal, W.~Lee, S.~W. Fung, and L.~Nurbekyan}, {\em Random features
  for high-dimensional nonlocal mean-field games}, Journal of Computational
  Physics, 459 (2022), p.~111136.

\bibitem{BaFe-NoDEA16}
{\sc M.~Bardi and E.~Feleqi}, {\em Nonlinear elliptic systems and mean-field
  games}, Nonlinear Differential Equations and Applications NoDEA, 23 (2016),
  pp.~1--32.

\bibitem{BaPr-SIAM14}
{\sc M.~Bardi and F.~S. Priuli}, {\em Linear-quadratic n-person and mean-field
  games with ergodic cost}, SIAM Journal on Control and Optimization, 52
  (2014), pp.~3022--3052.

\bibitem{ChFuLiNuOs-IP23}
{\sc Y.~T. Chow, S.~W. Fung, S.~Liu, L.~Nurbekyan, and S.~Osher}, {\em A
  numerical algorithm for inverse problem from partial boundary measurement
  arising from mean field game problem}, Inverse Problems, 39 (2023),
  p.~014001.

\bibitem{DiLiOsYi-JSC22}
{\sc L.~Ding, W.~Li, S.~Osher, and W.~Yin}, {\em A mean field game inverse
  problem}, Journal of Scientific Computing, 92 (2022), p.~7.

\bibitem{DiLiZh-arXiv23}
{\sc M.-H. Ding, H.~Liu, and G.-H. Zheng}, {\em Determining a stationary mean
  field game system from full/partial boundary measurement}, arXiv:2308.06688,
  (2023).

\bibitem{FuLiOsLi-JCP23}
{\sc G.~Fu, S.~Liu, S.~Osher, and W.~Li}, {\em High order computation of
  optimal transport, mean field planning, and potential mean field games}, J.
  Comput. Phys., 491 (2023), p.~112346.

\bibitem{GoPiSa-ESAIM16}
{\sc D.~A. Gomes, E.~Pimentel, and H.~S{\'a}nchez-Morgado}, {\em Time-dependent
  mean-field games in the superquadratic case}, ESAIM: Control, Optimisation
  and Calculus of Variations, 22 (2016), pp.~562--580.

\bibitem{ImLiYa-AML23}
{\sc O.~Imanuvilov, H.~Liu, and M.~Yamamoto}, {\em Unique continuation for a
  mean field game system}, Applied Mathematics Letters,  (2023), p.~108757.

\bibitem{ImLiYa-IPI24}
{\sc O.~Imanuvilov, H.~Liu, and M.~Yamamoto}, {\em Lipschitz stability for
  determination of states and inverse source problem for the mean field game
  equations}, Inverse Problems and Imaging,  (2024).

\bibitem{ImYa-arXiv23}
{\sc O.~Imanuvilov and M.~Yamamoto}, {\em Global {Lipschitz} stability for an
  inverse coefficient problem for a mean field game system}, arXiv:2307.04025,
  (2023).

\bibitem{KlAv-SIAM24}
{\sc M.~V. Klibanov and Y.~Averboukh}, {\em Lipschitz stability estimate and
  uniqueness in the retrospective analysis for the mean field games system via
  two {Carleman} estimates}, SIAM J. Math. Anal., 56 (2024), pp.~616--636.

\bibitem{KlLiLi-arXiv23C}
{\sc M.~V. Klibanov, J.~Li, and H.~Liu}, {\em Coefficient inverse problems for
  a generalized mean field games system with the final overdetermination},
  arXiv:2305.01065,  (2023).

\bibitem{KlLiLi-JIIP24}
\leavevmode\vrule height 2pt depth -1.6pt width 23pt, {\em On the mean field
  games system with lateral cauchy data via carleman estimates}, Journal of
  Inverse and Ill-posed Problems, 32 (2024), pp.~277--295.

\bibitem{KlLiYa-arXiv23}
{\sc M.~V. Klibanov, J.~Li, and Z.~Yang}, {\em Convexification numerical method
  for the retrospective problem of mean field games}, arXiv:2306.14404,
  (2023).

\bibitem{LaWo-TRB-2011}
{\sc A.~Lachapelle and M.-T. Wolfram}, {\em On a mean field game approach
  modeling congestion and aversion in pedestrian crowds}, Transportation
  Research Part B: Methodological, 45 (2011), pp.~1572--1589.

\bibitem{LeLiTeLiOs-SIAM21}
{\sc W.~Lee, S.~Liu, H.~Tembine, W.~Li, and S.~Osher}, {\em Controlling
  propagation of epidemics via mean-field control}, SIAM Journal on Applied
  Mathematics, 81 (2021), pp.~190--207.

\bibitem{LiMoZh-IP23}
{\sc H.~Liu, C.~Mou, and S.~Zhang}, {\em Inverse problems for mean field
  games}, Inverse Problems, 39 (2023), p.~085003.

\bibitem{LiZh-arXiv22}
{\sc H.~Liu and S.~Zhang}, {\em On an inverse boundary problem for mean field
  games}, arXiv preprint arXiv:2212.09110,  (2022).

\bibitem{LiZh-arXiv23}
\leavevmode\vrule height 2pt depth -1.6pt width 23pt, {\em Simultaneously
  recovering running cost and {Hamiltonian} in mean field games system}, arXiv
  preprint arXiv:2303.13096,  (2023).

\bibitem{LiZh-arXiv24}
\leavevmode\vrule height 2pt depth -1.6pt width 23pt, {\em Inverse boundary
  problem for a mean field game system with probability density constraint},
  arXiv:2402.13274,  (2024).

\bibitem{LiJaLiNuOs-SIAM21}
{\sc S.~Liu, M.~Jacobs, W.~Li, L.~Nurbekyan, and S.~J. Osher}, {\em
  Computational methods for first-order nonlocal mean field games with
  applications}, SIAM J. Numer. Anal., 59 (2021), pp.~2639--2668.

\bibitem{PiVo-IUMJ17}
{\sc E.~A. Pimentel and V.~Voskanyan}, {\em Regularity for second-order
  stationary mean-field games}, Indiana University Mathematics Journal,
  (2017), pp.~1--22.

\bibitem{ReSoWa-arXiv23}
{\sc K.~Ren, N.~Soedjak, and K.~Wang}, {\em Unique determination of cost
  functions in a multipopulation mean field game model}, Submitted,  (2023).
\newblock arXiv:2312.01622.

\bibitem{YuXiChLa-arXiv24}
{\sc J.~Yu, Q.~Xiao, T.~Chen, and R.~Lai}, {\em A bilevel optimization method
  for inverse mean-field games}, arXiv:2401.05539,  (2024).

\end{thebibliography}

}

\end{document}